\documentclass[12pt]{article}
\usepackage{array}
\usepackage{amsthm,amsmath,amssymb,color}
\usepackage{fullpage}

\usepackage{color}
\usepackage[normalem]{ulem}

\newtheorem{thm}{Theorem}[section]

\newtheorem{lem}[thm]{Lemma}
\newtheorem{prop}[thm]{Proposition}

\begin{document}

\title{ $K_{r+1}$-saturated graphs with small spectral radius}
\author{Jaehoon Kim\thanks{Mathematical Sciences Department, KAIST, jaehoon.kim@kaist.ac.kr}\, 
	Seog-Jin Kim\thanks{Department of Mathematics Education, Konkuk University, Seoul, 05029, Korea, skim12@konkuk.ac.kr} \thanks{This work was supported by the National Research Foundation of Korea(NRF) grant funded by the Korea government(MSIT)(NRF-2018R1C1B6003786).}\,
		Alexandr V. Kostochka\thanks{Department of Mathematics, University of Illinois, Urbana, IL, 61801, USA
		and
Sobolev Institute of Mathematics, Novosibirsk 630090, Russia, kostochk@math.uiuc.edu. Research of this author is supported in part by NSF grant
 DMS-1600592  and by grants 18-01-00353  and 19-01-00682 of the Russian Foundation for Basic Research.}\, and
Suil O\thanks{Department of Applied Mathematics and Statistics, The State University of New York, Korea, Incheon, 21985, suil.o@sunykorea.ac.kr. Research supported by NRF-2018K2A9A2A06020345 and by NRF-2020R1F1A1A01048226.}
}

\maketitle

\begin{abstract}

For a graph $H$, a graph $G$ is $H$-saturated if $G$ does not contain $H$ as a subgraph but for any 
$e \in E(\overline{G})$, $G+e$ contains $H$.
 In this note, we prove a sharp lower bound for the 
 number of paths and walks on length $2$ in $n$-vertex $K_{r+1}$-saturated graphs.
We then use this 
bound 
to give a lower bound on the spectral radii of such graphs which is asymptotically tight 
for each fixed $r$ and $n\to\infty$.

\bigskip


\end{abstract}

\section{Introduction}
\subsection{Notation and preliminaries}
In this note we  deal  with finite  undirected graphs with no loops or multiple edges.
For a graph $H$, a graph $G$ is $H$-{\em saturated}  if  $H$ is not a subgraph of
$G$ but after adding to $G$ any edge results in a graph containing $H$. 
For a positive integer $n$ and a graph $H$, the {\em extremal number} $ex(n, H)$ is the maximum number of edges in an $n$-vertex graph not containing $H$.
Clearly, an extremal $n$-vertex graph $G$ not containing $H$ with $|E(G)|=ex(n, H)$ is
$H$-saturated. Thus, one can also say that   $ex(n, H)$ is the maximum number of edges in an $n$-vertex  $H$-saturated graph. On the other hand,
the {\em saturation number of $H$},  $sat(n, H)$, is the least number of edges in an $H$-saturated graph with $n$ vertices.

Initiating the study of extremal graph theory, Tur{\'a}n~\cite{T1} determined the extremal number $ex(n, K_{r+1})$. He also proved that there is the unique extremal graph,  $T_{n,r}$, the $n$-vertex complete $r$-partite graph 
whose partite sets differ in size at most 1.
The first result on saturation numbers is due to Erd{\H{o}}s, Hajnal and Moon~\cite{EHM}:

\bigskip
{\bf Theorem A~\cite{EHM}.} {\em If $2\leq r<n$, then $sat(n, K_{r+1})=(r-1)(n-r+1)+\binom{r-1}{2}$. The only $n$-vertex $K_{r+1}$-saturated graph with
$sat(n, K_{r+1})$ edges is the graph $S_{n,r}$ obtained from
a copy of $K_{r-1}$ with vertex set $S$ by adding $n - r + 1$ vertices, each of which has neighborhood $S$. }

\medskip
Graph $S_{n,r}$ has clique number $r$ and no $r$-connected subgraphs; in particular, $S_{n,2}$ is a star. For an excellent survey
on saturation numbers, we refer the reader to Faudree, Faudree, and Schmitt~\cite{FFS}.

Recently, there was a series of publications on  eigenvalues of  $H$-free graphs.
For a graph $G$, let $A(G)$ be its adjacency matrix, and we index the eigenvalues of $A(G)$ in nonincreasing order, $\lambda_1(G)\ge \cdots \ge \lambda_n(G)$. 
 The value $\lambda_1(G)$ is also called the {\em spectral radius} of $G$, and denoted by $\rho(G)$.

Studying properties of quasi-random graphs, Chung, Graham, and Wilson~\cite{CGW}  proved a theorem implying that,
if $n$ is sufficiently large, $0 < c < \frac 12$ and $G$ is an $n$-vertex $K_r$-free graph with
 $\lceil cn^2 \rceil$ edges, then either $\lambda_n(G) < -c'n$ or $\lambda_2(G)>c'n$, where $c' = c'(r, c)$ is
a positive constant. However, the methods in~\cite{CGW} fail to indicate which of the two inequalities actually holds.
Bollob{\'a}s and Nikiforov~\cite{BN} observed that if $G$ is a dense $K_r$-free graph, then $\lambda_n(G) < - cn $ for some $c > 0$ independent of $n$.
Nikiforov~\cite{N1} gave a more precise statement that
if $G$ is a $K_{r+1}$-free graph with $n$ vertices and $m$ edges, then $\lambda_n(G) < -\frac{2^{r+1}m^r}{rn^{2r-1}}.$

Nikiforov~\cite{N2} also proved that if $G$ ia a $K_{r+1}$-free graph with $n$ vertices, then $\rho(G) \le \rho(T_{n,r}).$ Since each $K_{r+1}$-saturated graph is  $K_{r+1}$-free, his theorem implies the following.

\bigskip
{\bf Theorem B~\cite{N1}.} {\em
 If $G$ is a $K_{r+1}$-saturated graph with $n$ vertices, then  $$\rho(G) \le \rho(T_{n,r}).$$
}


In this note, we give  a new lower bound
for the spectral radius of an $n$-vertex $K_{r+1}$-saturated graph. This bound is asymptotically  tight 
when $r$ is fixed or grows as $o(n)$.
For this, we give a tight 
lower bound on the sum of the squares of the vertex degrees in an $n$-vertex $K_{r+1}$-saturated graph.


\subsection{Results}

Our main tool will be the following.

\begin{thm}
\label{main1}
 If $n\geq r+1$ and $G$ is a $K_{r+1}$-saturated graph with $n$ vertices,
 then  
 \begin{equation}\label{ne1}
 \sum_{v \in V(G)}d^2(v) \ge (n-1)^2(r-1)+(r-1)^2(n-r+1).
 \end{equation} 
For $r=2$, equality in the bound holds only when $G$ is $S_{n,2}$ or a Moore graph with diameter 2. 
For $r \ge 3$, equality in the bound holds only when $G$ is $S_{n,r}$.
\end{thm}

The reason why it is helpful is the following simple observation.

\begin{lem}\label{lem2} For every $n$-vertex graph $G$ with adjacency matrix $A$,
\begin{equation}\label{l2}
\rho^2(A)\geq \frac{1}{n}\sum_{v\in V(G)}d^2(v).
\end{equation}
\end{lem}

Theorem~\ref{main1} together with this observation immediately yield

\begin{thm}\label{main2}
If $2\leq r<n$ and $G$ is a $K_{r+1}$-saturated graph with $n$ vertices,
then 
\begin{equation}\label{so}
\rho(G) \ge \sqrt{\frac{(n-1)^2(r-1)+(r-1)^2(n-r+1)}{n}}.
\end{equation} 
\end{thm}

This bound asymptotically is tight 
because 
the spectral radius of $S_{n,r}$ is close to $f(n,r)$, where $f(n,r)$ is the lower bound for $\rho(G)$ in Theorem~\ref{main2}. More specifically, note that $\rho(S_{n,2})=f(2,n)$ and for $r \ge 3$, we have $\rho(S_{n,r})=f(r,n)+\frac{r-2}2+\Theta(\frac{r^{1.5}}{\sqrt{n}})$.

\begin{prop}\label{jk}
For  integers $2\leq r<n$,
$$\rho(S_{n,r})=\frac{r-2+\sqrt{(r-2)^2+4(r-1)(n-r+1)}}{2}.$$
\end{prop}

In the next section we prove Theorem~\ref{main1} (in a somewhat stronger form) and in the last section we present proofs for Lemma~\ref{lem2} and
Proposition~\ref{jk}

For undefined terms, see Brouwer and Haemers~\cite{BH},  Godsil and Royle~\cite{GR}, or West~\cite{W}.

\section{Proof of Theorem~\ref{main1}}

We will derive Theorem~\ref{main1} from the following slightly stronger statement.

\begin{thm}
\label{main3}
 If $n\geq r+1$ and $G$ is a $K_{r+1}$-saturated graph with $n$ vertices,
 then   \begin{equation}\label{ne3}
 \sum_{v \in V(G)}(d(v)+1)(d(v)+1-r) \ge (r-1)n(n-r). 
  \end{equation}
\end{thm}

\begin{proof}
	Let $m=|E(G)|$ and $\overline{m}=|E(\overline{G})|={n\choose 2}-m$. For $v\in V(G)$, let $f(v)$ be the number of pairs of non-adjacent vertices $x$ and $y$ in $N(v)$ such that $G[N(x)\cap N(y)\cap N(v)]$ contains 
	$K_{r-2}$ as a subgraph. Note that if $G[N(v)]$ is a copy of $K_{r-1}$, then $f(v)=0$.\\

\noindent
{\it Claim 1. $\overline{m}\leq \frac{1}{r-1}\sum_{v \in V(G)} f(v)$.}\\

 We construct an auxiliary bipartite graph $H$ with parts $A$ and $B$ as follows. Let $A=E(\overline{G})$  and  $B=V(G)$.
  The graph $H$ has an edge between $xy \in A$ and $v \in B$ iff $x,y \in N(v)$ and $G[N(x)\cap N(y)\cap N(v)]$ contains $K_{r-2}$ as a subgraph. 
   Then for each $v \in B$, we have $|N_{H}(v)|=f(v)$. Also, since $G$ is $K_{r+1}$-saturated, 
   for each $xy \in A$,  $G+xy$ contains $K_{r+1}$ as a subgraph. Thus there exist at least $r-1$ vertices $v$ such that $x, y \in N(v)$ and $G[N(x)\cap N(y)\cap N(v)]$ contains $K_{r-2}$ as a  subgraph, which implies 
\begin{equation}\label{eq2}
|N_{H}(xy)| \ge r-1.
\end{equation} 
By~\eqref{eq2},
\begin{equation}\label{eq21}
(r-1)\overline{m} \le \sum_{xy\in A} d_H(xy)= |E(H)|=\sum_{v \in V(G)}f(v).
 \end{equation} 
 This proves Claim 1. \\

\noindent
{\it Claim 2. For each $v \in V(G)$,  we have $\displaystyle f(v)\le {d(v)-r+2 \choose 2}$.}

Let $H_v=G[N(v)]$, and let $d(v)=p$. Since $G$ contains no $K_{r+1}$, the graph $H_v$ has no $K_{r}$. 
Partition the pairs of vertices in $N(v)$ into the sets $E_1,E_2$ and $E_3$ as follows:\\
(i) $E_1=E(H_v)$, \\
(ii) $E_2$ is the set  of the edges $xy\in E(\overline{H_v})$ such that $H_v+xy$ does not contain $K_r$,\\
(iii) $E_3$ is the set of the edges $xy\in E(\overline{H_v})$ such that $H_v+xy$ contains $K_r$. \\

Let $m_i=|E_i|$ for $1\leq i\leq 3$. By definition, $m_3=f(v)$ and 
 $m_1+m_2+m_3={p \choose 2}.$
As any $K_r$-free graph is a subgraph of $K_r$-saturated graph on the same vertex set, there exists a $K_r$-saturated graph $H'$ with vertex set $N(v)$ containing $H_v$. Then $E(H')\supseteq E_1$. Furthermore, since $H'$ is $K_{r}$-free and contains $E_1$, $E(H')\cap E_3=\emptyset$.
By Theorem A, $|E(H')|\geq (r-2)(p-r+2)+\binom{r-2}{2}$. Hence
$$m_3\leq {p\choose 2}-|E(H')|\leq  {p\choose 2}-(r-2)(p-r+2)-\binom{r-2}{2}={p-r+2 \choose 2}.$$
 This proves Claim 2. \\


Now we are ready to prove the theorem.
By Claims 1 and 2, 
$$
{n \choose 2}=m+\overline{m} \le m + \frac 1{r-1}\sum_{v\in V(G)} f(v) \le \sum_{v\in V(G)}\left[ \frac {d(v)}2+\frac 1{r-1}\frac{(d(v)-r+2)(d(v)-r+1)}2\right].
$$
Multiplying both sides by $2(r-1)$, we get
$$(r-1)n(n-1) \le \sum_{v \in V(G)}\left[(r-1)d(v)+(d(v)+1)(d(v)-r+1)-(r-1)(d(v)-r+1)\right].$$

This yields $$\sum_{v \in V(G)}(d(v)+1)(d(v)+1-r)\ge (r-1)n(n-1)-(r-1)^2n=(r-1)n(n-r),$$
and  Theorem~\ref{main3} is proved.
\end{proof}


To obtain  Theorem~\ref{main1}, observe that~\eqref{ne3} implies

$$\sum_{v \in V(G)}d^2(v)  \ge   (r-1)n(n-r)+(r-1)n+(r-2)2m.$$
So, by Theorem A,
\begin{equation} \label{eq3}
\sum_{v \in V(G)}d^2(v)  \ge (r-1)n(n-r)+(r-1)n+2(r-2)\left[{n \choose 2}-{n-r+1 \choose 2}\right] 
\end{equation}
$$ =  (n-1)^2(r-1)+(r-1)^2(n-r+1). $$
This proves the first part of Theorem~\ref{main1}.
Furthermore, for $r \ge 3$, equality in the bound requires equality in Theorem A.
Thus equality holds only for $S_{n,r}$.  

Suppose now $r=2$ and $G$ is an $n$-vertex $K_3$-saturated graph for which~\eqref{ne1} holds with equality. As $G$ is $K_3$-saturated, $G$ has diameter 2. 
Equality in the bound requires equality in~\eqref{eq21}, and hence equality in~(\ref{eq2}) for every $xy\in E(\overline{G})$. This means
$G$ has no $C_4$, which implies that $G$ has girth at least 5.
If $G$ has no cycles, then $G$ is a copy of $S_{n,2}$. Otherwise, $G$ is a Moore graph with diameter 2.

\medskip
Recall that there are at most four Moore graphs with diameter 2:  $C_5$, the Petersen graph, the Hoffman-Singleton graph, and possibly 
 one  57-regular graph of girth $5$ with 3250 vertices.

\section{Spectral radius}

We will use the following standard tool.

\begin{thm}[Rayleigh Quotient Theorem]\label{rqt}
	For a real matrix $A$ 
\begin{equation} \label{rqt1}	
	\rho(A)=\max_{x \in \mathrm{R}^n \setminus \{0\}}\frac{x^TAx}{x^Tx}.
	\end{equation}
\end{thm}

First, we present a proof of Lemma~\ref{lem2}. By~\eqref{rqt1},
$$\rho^2(A)=\rho(A^2)=\max_{x \in \mathrm{R}^n \setminus \{0\}}\frac{x^TA^2x}{x^Tx}=
\max_{x \in \mathrm{R}^n \setminus \{0\}}\frac{(x^TA^T)(Ax)}{x^Tx}\geq \frac{({\bf 1}^TA^T)(A{\bf 1})}{{\bf 1}^T{\bf 1}}=\frac{1}{n}\sum_{v\in V(G)}d^2(v).
$$
Thus,~\eqref{l2} holds. Together with Theorem~\ref{main1}, this implies Theorem~\ref{main2}.

\bigskip
To show that  Theorem~\ref{main2} is asymptotically tight, we will determine the spectral radius of $S_{n,r}$, i.e. prove
Proposition~\ref{jk}. We will need a new notion.
  Consider a partition $V(G) = V_1 \cup \cdots \cup V_s$ of the vertex set of a graph $G$ into $s$ non-empty subsets. For $1 \le i, j \le s$, let $q_{i,j}$ denote the average number
of neighbors in $V_j$ of the vertices in $V_i$. The quotient matrix $Q$ of this partition is the $s \times s$ matrix whose $(i, j)$-th entry equals $q_{i,j}$. The eigenvalues of the quotient matrix interlace the
eigenvalues of $G$. This partition is {\it equitable} if for each $1 \le i, j \le s$, each vertex $v \in V_i$ has exactly $q_{i,j}$ neighbors in $V_j$. In this case, the eigenvalues of the quotient matrix are eigenvalues of $G$ and the spectral radius of the quotient matrix equals the spectral radius of
$G$ (see \cite{BH}, \cite{GR} for more details).

\bigskip
\noindent [Proof of Proposition \ref{jk}] \\
Partition  $V(S_{n,r})$ into sets $A$ and $B$ such that $S_{n,r}[A]$ is a copy of $K_{r-1}$ and $S_{n,r}[B]$ is an independent set with $n-r+1$ vertices. Each vertex in $A$ is adjacent to all vertices in $B$. The quotient matrix of the partitions $A$ and $B$ is 
$$\begin{pmatrix}
r-2 & n-r+1 \\
r-1 & 0
\end{pmatrix}.$$
The characteristic polynomial of the matrix is $x^2 -(r-2)x -(r-1)(n-r+1)=0$. Since the partition $V(S_{n,r})=A \cup B$ is equitable,
$$\rho(S_{n,r})=\frac{r-2+\sqrt{(r-2)^2+4(r-1)(n-r+1)}}{2}.$$ 
This completes the proof of Proposition \ref{jk}.  \qed

\bigskip
Note that $\rho(S_{n,2})=\sqrt{n-1}$. Thus for $r=2$, equality in Theorem~\ref{main2} holds if and only if $G$ is $S_{n,2}$ or a Moore graph.
For $r \ge 3$, the bound in Theorem~\ref{main2} may be improved, and we guess that the spectral radius of $S_{n,r}$ is the minimum of $\rho(H)$ among all $n$-vertex $K_{r+1}$-saturated graphs $H$.

\bigskip
{\bf Acknowledgement.} We thank Xuding Zhu for helpful comments.

\end{document}